\documentstyle[twoside,12pt]{article}

\newcommand{\chapter}{\section}

\begin{document}

%Declaration section
%\theoremstyle{plain}
\newtheorem{Thm}{Theorem}%[section]
\newtheorem{Ax}{Axiom}%[section]
\newtheorem{Prop}{Proposition}%[subsection]
\newtheorem{Cor}[Prop]{Corollary}
\newtheorem{Main}{}
\renewcommand{\theMain}{}
\newtheorem{Lem}[Prop]{Lemma}
\newtheorem{Fact}{Fact}
\renewcommand{\theFact}{}

\newtheorem{Def}{Definition}%[section]
\newtheorem{rmk}{Remark}%[subsection]
\newenvironment{Rmk}{\begin{rmk}\em}{\end{rmk}}
\newtheorem{exm}{Example}%[subsection]
\newenvironment{Exm}{\begin{exm}\em}{\end{exm}}

\newcommand{\qed}{\par {\EM QED} }
\newtheorem{prf}{Proof}
\renewcommand{\theprf}{}
\newenvironment{Prf}{\begin{prf}\em}{\qed\end{prf}}
\newtheorem{prff}{}
\renewcommand{\theprff}{}
\newenvironment{Prff}{\begin{prff}\em}{\qed\end{prff}}
%\newtheorem{notation}{Notation}
%\renewcommand{\thenotation}{}

%Command section
%\errorcontextlines=0
%\numberwithin{equation}{section}
%\renewcommand{\rm}{\normalshape}%
% redefining \rm to mean: change to roman style

%\hsize=7.5in
%\vsize=10.5in

% Line spacing
%\renewcommand{\baselinestretch}{1.6}

% macroes

\newcommand{\YES}[1]{#1}
\newcommand{\NOT}[1]{}

\newcommand{\cA}{{\cal A}}
\newcommand{\cB}{{\cal B}}
\newcommand{\cC}{{\cal C}}
\newcommand{\cD}{{\cal D}}
\newcommand{\cE}{{\cal E}}
\newcommand{\cF}{{\cal F}}
\newcommand{\cG}{{\cal G}}
\newcommand{\cH}{{\cal H}}
\newcommand{\cI}{{\cal I}}
\newcommand{\cJ}{{\cal J}}
\newcommand{\cK}{{\cal K}}
\newcommand{\cL}{{\cal L}}
\newcommand{\cM}{{\cal M}}
\newcommand{\cN}{{\cal N}}
\newcommand{\cO}{{\cal O}}
\newcommand{\cP}{{\cal P}}
\newcommand{\cQ}{{\cal Q}}
\newcommand{\cR}{{\cal R}}
\newcommand{\cS}{{\cal S}}
\newcommand{\cT}{{\cal T}}
\newcommand{\cU}{{\cal U}}
\newcommand{\cV}{{\cal V}}
\newcommand{\cW}{{\cal W}}
\newcommand{\cX}{{\cal X}}
\newcommand{\cY}{{\cal Y}}
\newcommand{\cZ}{{\cal Z}}

\newcommand{\bbb}[1]{{\mbox{\bf #1}}}

\newcommand{\bN}{\bbb{N}}
\newcommand{\bZ}{\bbb{Z}}
\newcommand{\bR}{\bbb{R}}
\newcommand{\bC}{\bbb{C}}
\newcommand{\bQ}{\bbb{Q}}
\newcommand{\bT}{\bbb{T}}

\newcommand{\noind}[1]{{\setlength{\parindent}{0cm} #1}}
\newcommand{\parsk}{\par\medskip}

\newcommand{\varend}{\end{document}}

\newcommand{\LP}{\left(}  \newcommand{\RP}{\right)}
\newcommand{\LQ}{\left[}  \newcommand{\RQ}{\right]}
\newcommand{\LB}{\left\{} \newcommand{\RB}{\right\}}
\newcommand{\LA}{\left\langle} \newcommand{\RA}{\right\rangle}
\newcommand{\LS}{\left|}  \newcommand{\RS}{\right|}
\newcommand{\LN}{\left\|} \newcommand{\RN}{\right\|}
\newcommand{\bLP}{\bigl(}  \newcommand{\bRP}{\bigr)}

\newcommand{\ts}{{\otimes}} \newcommand{\TS}{{\bigotimes}}

\newcommand{\es}{\emptyset}
\newcommand{\sm}{\setminus} \newcommand{\st}{\subset}
\newcommand{\eps}{\varepsilon}

\newcommand{\beware}{~$\Psi\!\Psi\!\Psi$~}

\newcommand{\EM}{\em\bf}

\newcommand{\BE}{\begin{equation}}  \newcommand{\EE}{\end{equation}}
\newcommand{\BER}{$$\begin{array}}  \newcommand{\EER}{\end{array}$$}
\newcommand{\head}[1]{
             \par\bigskip\noind{{\large\bf#1}\nopagebreak\par}\medskip}
\newcommand{\chead}[1]{\begin{center}
              \par\bigskip{\large\bf#1}\nopagebreak\par\medskip
              \end{center}}

\newcommand{\fillitem}{\L{\embox{}}\vspace{-.6cm}}

\newcommand{\dfrac}[2]{\displaystyle{\frac{#1}{#2}}}
\newcommand{\toprove}[1]{{\buildrel\mbox{?}\over#1}}
\newcommand{\DEF}{{\buildrel\mbox{\scriptsize\,\,def\,\,}\over=}}
\newcommand{\BAR}[1]{\overline{#1}}
\newcommand{\I}{\infty}
\newcommand{\al}{\alpha} \newcommand{\be}{\beta}
\newcommand{\OM}{\Omega} \newcommand{\om}{\omega}
\newcommand{\la}{\lambda}
\newcommand{\mod}{\mbox{ mod }}
\newcommand{\half}{\frac{1}{2}}
\newcommand{\NI}{{n\to\infty}}
\newcommand{\ang}{{<\!\!\!)}}
\newcommand{\arc}[1]{\breve{#1}}

\newenvironment{txeqn}[1]{\begin{equation}\label{#1}\begin{array}{l}}{
                         \end{array}\end{equation}}

%\newcommand{\Th}[2]{\par\medskip{\EM#1}. {\em#2}\par\medskip}

%. \newpage . \newpage
%\setcounter{page}{1}

%special newcommands for this article
\newcommand{\supp}{{\mbox{supp}\,}}
\newcommand{\spc}{{\mbox{Spec}\,}}
\newcommand{\cs}{continuous}
\newcommand{\fn}{function}
\newcommand{\Bh}{Banach}
\newcommand{\HAT}[1]{\widehat{#1}}

%Topmatter

\title{Differentiability of Scalar Functions Applied to Hermitian
Operators -- a Fourier Transform Approach}

%Author info
\author{Eliahu Levy\\
%\address{
Department of Mathematics\\
Technion -- Israel Institute of Technology,
Haifa 32000, Israel\\
email: eliahu@techunix.technion.ac.il}

%\thanks{}
%\email{@}
%\keywords{}
%\subjclass{Primary 46G05, 47B15}
\date{}
%End topmatter

\maketitle

\begin{abstract}
Let $g$ be a (say, sufficiently differentiable) scalar \fn\ on
the reals. One knows how to apply $g$ to Hermitian elements $A$
of a $C^*$-algebra\NOT{(Everything generalizes easily to
multivariate \fn s applied to commuting normal elements)}. Yet the
question of differentiability of $A\mapsto g(A)$ is not trivial,
since in general ``$A$ and $dA$ do not commute''. However, since
the mapping from $g$ to $A\mapsto g(A)$ is linear, one can,
via Fourier Transform, reduce the case of general $g$ to the
case $g(x)=\exp(x)$. For the latter one has an explicit formula
for the $n$-th Fr\'echet derivative (more complicated than in the
scalar case -- still $A$ and $dA$ do not commute!). In this way,
one bounds the norm of the $n$-th (Fr\'echet) derivative of
$A\mapsto g(A)$ on a ball of radius $r$ by a Sobolev norm involving the
$(n+1)$-th derivative of $g$ on the interval $[-r;r]$.
\end{abstract}

\tableofcontents

\section{Derivatives of Functions Between Banach Spaces}
\label{Sec1}
We assume known the notion of integral of a \fn\ having values in
a locally convex space $E$ ($\int f$ being defined as the member
of $E$ such that for any \cs\ linear \fn al $\tilde{x}$ on $E$
$\tilde{x}(\int f)=\int\tilde{x}\circ f$), as well as the theorems
ensuring the existence of the integral if $E$ is quasi-complete
(in particular, a Fr\'echet space),
and either $f$ is \cs\ from a compact space $K$ to $E$ and the integration is
on a finite measure on $K$ (i.e.\ with finite mass), or $f$ is
continuous from a locally compact space to $E$ and the integration
is on a positive measure, finite on any compact and such that for
each topology-defining seminorm $\rho$ on $E$ $\int\rho\circ
f<\I$.
\parsk

For $E$ a \Bh\ space, $x\in E$ and $r>0$ denote
$$B_r:=\{y\in E\,|\,\|y\|<r\},\qquad B_r(x):=\{y\in
E\,|\,\|y-x\|<r\}$$ If $E$ and $F$ are \Bh\ spaces denote by
$\cL(E,F)$ the \Bh\ space of the bounded linear operators from $E$
to $F$. For $n\ge0$ integer denote by $\cL^{(n)}(E,F)$ the \Bh\
space of bounded $n$-multi-linear operators from $E^n$ to $F$. One
may identify $\cL^{(0)}(E,F)=F$, $\cL^{(1)}(E,F)=\cL(E,F)$,
$\cL^{(n+1)}(E,F)=\cL\LP E,L^{(n)}(E,F)\RP$.
\parsk

We shall need the following
\parsk

\begin{Prop}
\label{Prop1.1} Let $E$ and $F$ be \Bh\ spaces, let $U\subset E$
be open and let $f:U\to F$, $g:U\to \cL(E,F)$ be \cs. Then a
necessary and sufficient condition that $g$ is the Fr\'echet
derivative of $f$ on $U$ is that for any $x\in U$ and for any
$v\in E$ in some open ball around $x$ contained in $U$ one has
\BE f(x+v)-f(x)=\int_0^1<g(x+tv),v>\,dt\label{(1.1)}.\EE
\end{Prop}

\begin{Prf}
If $g$ is even the G\^ateau derivative of $f$ on $U$,
$x\in U$ and $\eps>0$ such that $B_\eps(x)\subset U$, then for all
$v\in B_\eps(x)$ the function $t\mapsto f(x+tv)$ is \cs\ on
$[0,1]$ and has on $]0,1[$ the \cs\ derivative
$t\mapsto<g(x+tv),v>$. Hence $(\ref{(1.1)})$ follows from the fundamental
theorem of calculus (for \Bh-space valued \fn s).
\parsk

Conversely, if $(\ref{(1.1)})$ holds for all $x\in U$ for $v$ in some open
ball $B_\eps(x)$ around $x$, then one has
$$\|f(x+v)-f(x)-<g(x),v>\|=\|\int_0^1<g(x+tv)-g(x),v>\,dt\|\le
\LP\sup_{\|y-x\|\le\|v\|}\|g(y)-g(x)\|\RP\|v\|,$$
and by the continuity of $g$ at $x$ and the definition of Fr\'echet
derivative one has $f'(x)=g(x)$.
\end{Prf}
\parsk

For $K$ a (Hausdorff) compact topological space and $F$ a \Bh\
space denote by $\cC(K,F)$ the \Bh\ space of continuous \fn s
$K\to F$ with the sup norm.
\parsk

\begin{Prop}
\label{Prop1.2} Let $K$ be Hausdorff compact and let $E$, $F$ be
\Bh\ spaces. Let $U\subset E$ be open and $f:U\times K\to F$ \cs.
Then the function $\tilde{f}: U\to\cC(K,F)$ given by
$\tilde{f}(x):=(s\mapsto f(x,s))$ is \cs.
\end{Prop}

\begin{Prf}
Let $x\in U$ and $\eps>0$. For each $s\in K$ there is
an open neighborhood $W_s$ and a $\delta_s>0$ so that $s'\in W_s$
and $\|x'-x\|<\delta_s$ imply $\|f(x',s')-f(x,s)\|<\eps$, hence,
since also $\|f(x,s')-f(x,s)\|<\eps$ one gets
$\|f(x',s')-f(x,s')\|<2\eps$. Since $K$ is compact, there exist
$s_1,\ldots,s_n\in K$ with $\cup W_{s_i}=K$. Take
$\delta=\min_i\delta_{s_i}$. Then for $\|x'-x\|<\delta$ we have
for all $s'\in K$ $\|f(x',s')-f(x,s')\|<2\eps$ thus
$$\|\tilde{f}(x')-\tilde{f}(x)\|=\sup_{s\in
K}\|f(x',s)-f(x,s)\|\le 2\eps.$$
\end{Prf}
\parsk

\begin{Prop}
\label{Prop1.3} Retaining the setting of Prop.\ \ref{Prop1.2}, assume $f$
has a (``partial'') Fr\'echet derivative w.r.t.\ $x$ everywhere in
$U\times K$, and this derivative, as a function $g:U\times K\to
L(E,F)$ is \cs. Then $\tilde{f}$ has a \cs\ derivative in $U$,
that derivative being the \fn\ $\tilde{g}:U\to L(E,\cC(K,F)$ given
by $\tilde{g}(x):=v\mapsto(s\mapsto<g(x,s),v>)$ ($v\in E$, $s\in
K$).
\end{Prop}

\begin{Prf}
Firstly, for each $v\in E$ the \fn\
$s\mapsto<g(x,s),v>$ is \cs\ on $K$ hence belongs to $\cC(K,F)$.
Moreover, $\tilde{g}$ is \cs\ from $U$ to $L(E,\cC(K,F)$ since
\begin{eqnarray*}
&&\|<\tilde{g}(x'),v>(s)-<\tilde{g}(x),v>(s)\|=
\|<g(x',s),v>-<g(x,s),v>\|\le\\
&&\LP\sup_{s\in K}\|g(x',s)-g(x,s)\|\RP\|v\|
\end{eqnarray*}
and the last supremum tends to $0$ as $x'\to x$ by Prop.\ \ref{Prop1.2}.
\parsk

Thus, to prove $D\tilde{f}=\tilde{g}$ we need only to prove the
equality of Prop.\ \ref{Prop1.1}:
$$f(x+v)-f(x)=\int_0^1<g(x+tv),v>\,dt,\eqno{(\ref{(1.1)})}$$
which will hold if substituting any $s$ would give equality, i.e.\
if
$$f(x+v,s)-f(x,s)=\int_0^1<g(x+tv,s),v>\,dt$$
which follows from Prop.\ \ref{Prop1.1}, $g$ being the ``partial'' derivative
of $f$.
\end{Prf}
\parsk

\begin{Cor}
\label{Cor1.4} If, in the setting of Prop.\ \ref{Prop1.2} and
\ref{Prop1.3}, $\mu$ is a finite measure on $K$ then in $U$:
$$<\dfrac{d}{dx}\int f(x,s)\,d\mu(s),v>=\int<g(x,s),v>\,d\mu(s).$$
\end{Cor}

\begin{Prf}
This follows from Prop.\ \ref{Prop1.3} and from the fact that if
one denotes, for $h\in\cC(K,F)$,\, $\mu(h):=\int h(s)\,d\mu(s)$,
one has
\begin{eqnarray*}
&&\int f(x,s)\,d\mu(s)=\mu(\tilde{f}(x))\\
&&\int<g(x,s),v>\,d\mu(s)=\int<\tilde{g}(x),v>(s)\,d\mu(s)=
\mu(<\tilde{g}(x),v>).
\end{eqnarray*}
\end{Prf}
\parsk

Let now $E$ and $F$ be \Bh\ spaces, and denote by $\cB\cC^\I(E,F)$
the linear space of \fn s $f:E\to F$ with \cs\ derivatives of any
order on all of $E$, these derivatives bounded on any bounded
subset of $E$, with the locally convex topology defined by the
seminorms
\BE\rho_{n,r}(f):=\sup_{\|x\|<r}\|f^{(n)}(x)\|=
\sup_{\|x\|\le r}\|f^{(n)}(x)\|\qquad r>0,\; n=0,1,2,\ldots\label{(1.2)}\EE
(the two suprema are equal since $f^{(n)}$ is \cs\ on $E$).
\parsk

This space is metrizable, being Hausdorff with topology defined by
a countable set of seminorms (take also $r>0$ integer). Also, this
space is topologically isomorphic to the subspace $Y$ of the
product $X$ of the \Bh\ spaces $X_{n,r}$ ($n\ge0$ and $r>0$
integers), $X_{n,r}$ being the space of \cs\ bounded \fn s from
$B_r$ to $\cL^{(n)}(E,F)$ with sup norm, where an element of the
product of $X_{n,r}$ is defined to be in $Y$ if it satisfies
compatibility conditions: that for any $r$ and $n$ the $n+1,r$
coordinate be the derivative of the $n,r$ coordinate and the $n,r$
coordinate coincide with the restriction of the $n,r+1$ coordinate
to $B_r$. All these conditions give closed sets (by Prop.\ \ref{Prop1.1}).
Hence $Y$ is closed in the product $X$, so $\cB\cC^\I(E,F)$ is
complete.
\parsk

One concludes that $\cB\cC^\I(E,F)$ is a Fr\'echet space.
\parsk

For $E=F$ denote $\cB\cC^\I(E):=\cB\cC^\I(E,E)$.

\section[Polynomials and Entire Functions]{Derivatives of Polynomials and
Entire Functions in a \Bh\ Algebra}
\label{Sec2}
From now on, $A$ will be a \Bh\ algebra with unit (denoted by
$1$), we also assume $\|1\|=1$,
$\|xy\|\le\|x\|\|y\|$,\,\,$x,y\in A$.
\parsk

Denote by $S_n$ ($n\ge0$ integer) the set of permutations in
$\{1,\ldots,n\}$. $S_n$ has $n!$ elements.
\parsk

Denote by $P_{n,k}$ ($n\ge0$, $k$ integers) the set of sequences
$\al=(\al_0,\ldots,\al_n)$ with $\al_i\ge0$ integers and
$\sum\al_i=k$.
\parsk

Let us find how many members $P_{n,k}$ has. $P_{n,k}$ has 1-1
correspondence with the set of sequences
$\beta=(\beta_1,\ldots,\beta_n)$ with $\beta_i$ integers,
$\beta_i\le\beta_{i+1}$, $\beta_1\ge0$, $\beta_n\le k$ (just
correspond $\al_i=\beta_{i+1}-\beta_i$ for $1\le i\le n-1$,
$\al_0=\beta_0$, $\al_n=k-\beta_n$.) The set of $\beta$'s has 1-1
correspondence with the set of $\gamma=(\gamma_1,\ldots,\gamma_n)$
with $\gamma_i$ integers, $\gamma_i<\gamma_{i+1}$, $\gamma_1\ge1$,
$\gamma_n\le n+k$ (just let $\gamma_i=\beta_i+i$). The $\gamma$'s
can be identified with the subsets of $n$ elements in
$\{1,\ldots,n+k\}$. One concludes that $P_{n,k}$ has ${n+k\choose
n}$ elements.
\parsk

\begin{Thm}
\label{Th2.1} For every integer $k\ge0$ the \fn\ $x\mapsto
x^k$ on $A$ is $\cC^\I$ and its $n$-th derivative is given by
($v_i\in A$):
\BE<D^n(x^k), v_1\otimes v_2\otimes\cdots\otimes v_n>=
\sum_{\phi\in S_n}\sum_{\al\in P_{n,k-n}}
x^{\al_0}v_{\phi(1)}x^{\al_1}v_{\phi(2)}\cdots
x^{\al_{n-1}}v_{\phi(n)}x^{\al_n}\label{(2.1)}\EE
(if $k-n<0$\, $P_{n,k-n}$ is empty and the sum is $0$).
\end{Thm}

\begin{Prf}
The function $x\mapsto x^k$ is clearly \cs\ (and
bounded on bounded subsets of $A$).
\parsk

Let $k>0$. To obtain the first derivative of $x\mapsto x^k$, note
that it is a composition of the multilinear \fn\ from $A^k$ to
$A$\, $(x_1,x_2,\ldots,x_n)\mapsto x_1x_2\cdots x_n$ and the
linear $x\mapsto(x,x,\ldots x)$. Hence $x^k$ is $\cC^1$ with the
\cs\ derivative given by
$$<D(x^k),v>=\sum_{i=1}^k x^{i-1}vx^{k-i}=\sum_{\al\in
P_{1,k-1}}x^{\al_0}vx^{\al_1}$$
agreeing with $(2.1)$ for $n=1$.
\parsk

For $k=0$ and $n>1$ we have $x^0=1$ so $D^n(x^0)=0$, and the RHS
in $(2.1)$ is also $0$ since $P_{n,-n}=\es$.
\parsk

The assertion holds also for $n=0$. Indeed, then the LHS is $x^k$,
and in the RHS $S_0=\es$, $P_{0,k}=\{(k)\}$ so the RHS is also
$x^k$.
\parsk

Proceed now by induction, assuming that $n\ge1$, $x^k$ is $\cC^n$
and $(\ref{(2.1)})$ holds for $n$. Firstly $(\ref{(2.1)})$ for $n$
implies that
$D^n(x^k)$ obtains from composition of a bounded multilinear \fn\
with the $\cC^1$ \fn s $x^{\al_i}$, so $D^n(x^k)$ is $\cC^1$,
making $x^k$ $\cC^{n+1}$. To prove $(\ref{(2.1)})$ for $n+1$,
differentiate the expression in the RHS in $(\ref{(2.1)})$ for $n$, using
the same $(\ref{(2.1)})$ for $n=1$ which we obtained above:
\begin{eqnarray*}
&&<D^{n+1}(x^k),v_1\otimes\cdots\otimes v_n\otimes w>=\\
&&\sum_{\phi\in S_n}\sum_{\al\in P_{n,k-n}}\sum_{i=1}^n
x^{\al_0}v_{\phi(1)}\cdots x^{\al_{i-1}}v_{\phi(i)}
\LP\sum_{j=1}^{\al_i}x^{j-1}wx^{\al_i-j}\RP
v_{\phi(i+1)}x^{a_{i+1}}\cdots v_{\phi(n)}x^{\al_n}.
\end{eqnarray*}
And one easily convinces oneself that the RHS here is exactly the
RHS of $(\ref{(2.1)})$ for $n+1$, with $v_{n+1}=w$, where the summand
corresponding to a $\phi\in S_n$, an $0\le i\le n$, an $\al\in
P_{n,k-n}$ and a $1\le j\le\al_i$ would correspond in the RHS of
$(\ref{(2.1)})$ for $n+1$ to the $\psi\in S_{n+1}$ with
$$\psi(m)=\left\{\begin{array}{ll}
\phi(m)&1\le m\le i\\
n+1&m=i+1\\
\phi(m-1)&i+1<m\le n+1
\end{array}\right.$$
and to the sequence
$$(\al_0,\ldots,\al_{i-1},j-1,\al_i-j,\al_{i+1},\ldots,\al_n)\in
P_{n+1,k-n-1}.$$
\end{Prf}
\parsk

\begin{Cor}
\label{Cor2.2} For integer $k\ge0$ the \fn\ $x\mapsto x^k$
is in $\cB\cC^\I(A)$ with seminorms (for $m<0$ define
$\dfrac1{m!}=0$):
\BE\rho_{n,r}(x^k)\le\dfrac{k!}{(k-n)!}\,r^{k-n}.\label{(2.2)}\EE
\end{Cor}

\begin{Prf}
This follows from Thm.\ \ref{Th2.1}, from $S_n$ having $n!$
elements and $P_{n,k-n}$\, ${k\choose n}$ elements, and from the
inequality, holding for $\|x\|\le r$ (with the notation of
$(\ref{(2.1)})$):
$$\LN x^{\al_0}v_{\phi(1)}x^{\al_1}v_{\phi(2)}\cdots
x^{\al_{n-1}}v_{\phi(n)}x^{\al_n})\RN\le
r^{k-n}\|v_1\|\|v_2\|\cdots\|v_n\|$$
\end{Prf}
\parsk

\begin{Thm}
\label{Th2.3} Let $(a_k)_{k=0,1,\ldots}$ be a sequence of complex
numbers with $|a_k|^{1/k}\to0$ (i.e.\ so that the complex \fn\
$\sum a_kz^k$ is entire). Then the series for $x\in A$
\BE f(x)=\sum_{k=0}^\I a_kx^k \label{(2.3)}\EE
converges absolutely for all $x\in A$, converges in $\cB\cC^\I(A)$
as a series of elements of that space hence defines an element
(\fn) belonging to $\cB\cC^\I(A)$.
\end{Thm}

\begin{Prf}
$\cB\cC^\I(A)$ being a Fr\'echet space, to prove
convergence of the series in $\cB\cC^\I(A)$ it suffices to prove
``absolute convergence'' w.r.t.\ the defining seminorms, i.e.\ to
prove that for $n\ge0$ integer and $r>0$ the following series
converges:
$$\sum_k\,\rho_{n,r}(a_kx^k)=\sum_k|a_k|\,\rho_{n,r}(x^k)\le^{(\ref{(2.2)})}
\sum_k\,|a_k|\dfrac{k!}{(k-n)!}r^{k-n}\le
r^{-n}\sum_k|a_k|\,k^nr^k.$$
And the last series converges since
$$\LP|a_k|\,k^nr^k\RP^{1/k}=|a_k|^{1/k}(k^{1/k})^nr\to_{k\to\I}0.$$
The absolute convergence of the series $(\ref{(2.3)})$ at each $x\in A$
(to the value at $x$ of the element that the series defines in
$\cB\cC^\I(A)$) follows from the fact that evaluation at $x$ is a
\cs\ operator on $\cB\cC^\I(A)$ (because for $g\in\cB\cC^\I(A)$, if
$\|x\|<r$ then $\|g(x)\|\le\rho_{0,r}(g)$).
\end{Prf}
\parsk

\begin{Rmk}
\label{Rmk2.4} Retaining the setting of Thm.\ \ref{Th2.3}, if one
replaces $a_k$ by $a_kt^k$ for some complex $t$, then also
$|a_kt^k|^{1/k}=|t||a_k|^{1/k}\to_{k\to\I}0$, so the \fn\
\BE x\mapsto\sum_{k=0}^\I a_kt^kx^k=f(tx)\label{(2.4)}\EE
is in $\cB\cC^\I(A)$. Moreover, the computations in seminorms in
the proof of Thm.\ \ref{Th2.3} show that the series $(\ref{(2.4)})$ converges,
as a series of elements of $\cB\cC^\I(A)$, uniformly in $t$ on every
disk $|t|\le R$. Hence the $\cB\cC^\I(A)$-values \fn\ on the
complex plane $t\mapsto(x\mapsto f(tx))$ is \cs.
\end{Rmk}

\section{The Derivatives of the Exponential Function}
\label{Sec3}
Take in Thm.\ \ref{Th2.3} $a_k=\dfrac1{k!}$. For any fixed natural $m$, for
$k>m$
$$|a_k|^{1/k}=\LP\dfrac1{k!}\RP^{1/k}\le\LP\dfrac1{m!k^{k-m}}\RP^{1/k}\le
\LP\dfrac1{m!}\RP^{1/k}\LP k^{1/k}\RP^m\dfrac1{k}\to_{k\to\I}0$$
hence $|a_k|^{1/k}\to_k0$ -- of course the complex \fn\
$\exp z=\sum\dfrac{z^k}{k!}$ is entire. By Thm.\ \ref{Th2.3}, the same
series defines a \fn\ $A\to A$, also denoted by $\exp$:
\BE\exp x:=\sum_{k=0}^\I\dfrac{x^k}{k!}\label{(3.1)},\EE
and this \fn\ is a member of $\cB\cC^\I(A)$.
\parsk

Clearly $\exp 0=1$ and for $x\in A$ any element of $x$ that
commutes with $x$ commutes with $\exp x$. For commuting $x,y\in A$
consider the family of elements of $A$
$\LP\dfrac1{k!m!}x^ky^m\RP_{k,m=0,1,\ldots}$. The sum of norms of
this family $<\I$, hence this family is unconditionally summable
in $A$. Grouping its elements varying $k$ then varying $m$ one
finds that its sum is $(\exp x)(\exp y)$, while by grouping
$(k,m)$ according to $k+m$ one gets $\exp(x+y)$ for the sum. Thus
the exponent in $A$ has the characteristic property, for {\em
commuting} $x,y\in A$:
\BE\exp(x+y)=(\exp x)(\exp y)\label{(3.2)}.\EE
\parsk

For $x\in A$ and complex $t$ one has
$$\exp(tx)=\sum_{k=0}^\I\dfrac{(tx)^k}{k!}=\sum_{k=0}^\I\dfrac1{k!}x^kt^k.$$
We wish to take derivatives w.r.t.\ $t$. The series of
$t$-derivatives of the terms gives:
$$\sum_{k=1}^\I\dfrac1{k!}x^kkt^{k-1}=\sum_{k=1}^\I\dfrac1{(k-1)!}x^kt^{k-1}=
x\sum_{k=0}^\I\dfrac1{k!}x^kt^k.$$
This series sums to $x(\exp(tx))=(\exp(tx))x$. Since each term and the sum
are \cs\ in $t$ for both the series for $\exp(tx)$ and the series
of term-by-term derivatives, and both series converge uniformly on
closed disks in $t$ (by Remark \ref{Rmk2.4}), and since the relation of
being the derivative behaves for \cs\ \fn s as in Prop.\ \ref{Prop1.1}, one
finally concludes:
\BE\dfrac{d}{dt}\,\exp(tx)=x(\exp(tx))=(\exp(tx))x.\label{(3.3)}\EE
\parsk

We wish to compute the higher derivatives of $\exp x$ w.r.t.\ $x$.
Contrary to the complex case, the problem here, as mentioned
above, is that ``$x$ and $dx$ do not commute''. One might try to use
the series $(\ref{(3.1)})$, but we shall proceed differently.
\parsk

\begin{Prop}
\label{Prop3.1} For $x,v\in A$:
$$<D(\exp x),v>=\int_0^1\exp(tx)v\exp((1-t)x)\,dt.$$
\end{Prop}

\begin{Prf}
Consider the $A$-valued \fn\ of two real variables:
$$f(s,t)=\exp(t(x+sv))\exp(-tx)$$
We know that $f$ is $\cC^\I$. One has, by $(\ref{(3.3)})$:
\begin{eqnarray*}
&&\dfrac{\partial^2}{\partial t\,\partial s}f(s,t)=
\dfrac\partial{\partial s}\dfrac\partial{\partial t}f(s.t)=\\
&&\dfrac\partial{\partial s}\LQ\exp(t(x+sv))(x+sv)\exp(-tx)+
\exp(t(x+sv))(-x)\exp(-tx)\RQ=\\
&&\dfrac\partial{\partial s}\LQ s\exp(t(x+sv))v\exp(-tx)\RQ=\\
&& s\dfrac\partial{\partial s}\LQ\exp(t(x+sv))\RQ v\exp(-tx)+
\exp(t(x+sv)v\exp(-tx).
\end{eqnarray*}
Substitute $s=0$:
\BE\dfrac\partial{\partial t}\LP\dfrac\partial{\partial s}
f(s,t)|_{s=0}\RP=\exp(tx)v\exp(-tx).\label{(3.4)}\EE
Also $\dfrac\partial{\partial s}f(s,0)|_{s=0}=
\dfrac\partial{\partial s}1=0$, hence $(\ref{(3.4)})$ gives:
\BE\dfrac\partial{\partial s}f(s,1)|_{s=0}=
\int\exp(tx)v\exp(-tx)\,dt\label{(3.5)}.\EE
But by the definition of $f$:
$$\dfrac\partial{\partial s}f(s,1)|_{s=0}=
\LP\dfrac\partial{\partial s}\exp(x+sv)|_{s=0}\RP \exp(-x)=
<D(\exp x),v>\exp(-x),$$
which together with $(\ref{(3.5)})$ gives our assertion.
\end{Prf}
\parsk

\begin{Cor}
\label{Cor3.2} For $s\ge 0$, $x,v\in A$:
$$<\dfrac{d}{dx}(\exp(sx)),v>=\int_0^s\exp(tx)v\exp((s-t)x)\,dt.$$
\end{Cor}

\begin{Prf}
\begin{eqnarray*}
&&<\dfrac{d}{dx}(\exp(sx)),v>=<\dfrac{d}{dy}(\exp y)|_{y=sx},sv>=\\
&&=s\int_0^1\exp(t_1sx)v\exp((1-t_1)sx)\,dt_1=_{t=st_1}
\int_0^s\exp(tx)v\exp((s-t)x)\,dt
\end{eqnarray*}
\end{Prf}
\parsk

We shall need the {\EM $n$-simplex} $\Delta_n$. Let
$\bR^{n+1}:=\bR^{\{0,1,\ldots,n\}}$.
Let $L_n$ be the $n$-dimensional hyperplane
$$L_n:=\{t\in\bR^{n+1}\,|\,\sum_{i=0}^nt_i=0\}.$$
Let $\omega_n$ be a translation-invariant $n$-form in $\bR^{n+1}$ such that
\BE\omega_n\wedge\LP\sum_{i=0}^n\,dt_i\RP=\pm\bigwedge_{i=0}^ndt_i.
\label{(3.6)}\EE
This does not determine $\omega_n$ uniquely, but determines the Haar
measure that $\omega_n$ induces on hyperplanes parallel to $L_n$. We
call it {\EM Lebesgue measure} on such hyperplanes and denote it by
$\lambda_n$.
\parsk

The $n$-simplex $\Delta_n$ is defined as:
$$\Delta_n:=\{t\in\bR^{n+1}\,|\,\forall_i\,t_i\ge0,\,\sum_{i=0}^nt_i=1\}.$$
It lies on a hyperplane parallel to $L_n$.
\parsk

For each $n\ge1$,\,$1\le k\le n$ we have a linear mapping
$T_{n,k}:\bR^{n+1}\to\bR^n$:
$$T_{n,k}(t_0,t_1,\ldots,t_n)=(t_0,\ldots,t_{k-2},t_{k-1}+t_k,t_{k+1},\ldots
t_n).$$
$T_{n,k}$ maps $\Delta_n$ onto $\Delta_{n-1}$. If we denote by
$T_{n,k}^*\omega_{n-1}$ the pull-back of $\omega_{n-1}$ (which is a
translation-invariant $n-1$ -form in $\bR^{n+1}$), then by pulling back
$(\ref{(3.6)})$ one gets:
$$T_{n,k}^*\omega_{n-1}\wedge\LP\sum_{i=0}^n\,dt_i\RP=
\pm(dt_{k-1}+dt_k)\wedge\bigwedge_{i\ne k-1,k} dt_i.$$
Hence if $\eta_{n,k}$ is any translation-invariant $1$-form on $\bR^{n+1}$
such that 
\BE\eta_{n,k}\wedge(dt_{k-1}+dt_k)=\pm dt_{k-1}\wedge dt_k\label{(eta)}\EE
(say, $\eta_{n,k}=dt_{k-1}$ or $\eta_{n,k}=dt_k$),
then $T_{n,k}^*\omega_{n-1}\wedge\eta_{n,k}$ may serve as $\omega_n$.
Thus

\begin{Cor} \label{Cor:eta}
Let $\eta_{k,n}$ as above satisfy (\ref{(eta)}). Then to integrate over
$\Delta_n$ w.r.t.\ $\lambda_n$, one may integrate over the fibers of
$T_{n,k}$ w.r.t.\ $|\eta_{n,k}|$ and then integrate over $\Delta_{n-1}$
w.r.t.\ $\lambda_{n-1}$.
\end{Cor}
\parsk

\begin{Prop}
\label{Prop3.3}
$$\lambda_n(\Delta_n)=\dfrac1{n!}.$$
\end{Prop}

\begin{Prf}
for $n=0$ $\Delta_0=\{1\}$ and $\omega_0=1$ (note that $\omega_0$, being a
translation-invariant $0$-form, may be viewed as a real number). $L_0$ is
a singleton, to which $\lambda_0$ gives the mass $1$, hence
$\lambda_0(\Delta_0)=1$.
\parsk

To pass from $n-1$ to $n$ ($n\ge1$), take $\eta_{n,n}=dt_n$ to have
$$\omega_n=T_{n,n}^*\omega_{n-1}\wedge dt_n.$$
Hence to integrate $\lambda_n$ over $\Delta_n$, one may integrate
$T_{n,n}^*\omega_{n-1}$ over the slices $t_n=1-s$ and then integrate
over $|ds|$ on $s\in[0,1]$. These slices are
$$\{t\in\bR^{n+1}\,|\,(t_0,\ldots,t_{n-1})\in s\cdot\Delta_{n-1},\,t_n=s-1\}$$ 
and there, for $t\in s\cdot\Delta_{n-1}$,
$$T_{n,n}(t_0,\ldots,t_{n-1},1-s)=(t_0,\ldots,t_{n-1}+1-s)$$
so integration w.r.t.\ $T_{n,n}^*\omega_{n-1}$ transfers to integration
w.r.t.\ $\lambda_{n-1}$ over $s\cdot\Delta_{n-1}$. Hence one concludes:
$$\lambda_n(\Delta_n)=\int_0^1\lambda_{n-1}(s\cdot\Delta_{n-1})\,ds=
\lambda_{n-1}(\Delta_{n-1})\cdot\int_0^1s^{n-1}\,ds=
\dfrac1n\lambda_{n-1}(\Delta_{n-1})$$
and our assertion follows.
\end{Prf}
\parsk

Combining Corollaries \ref{Cor3.2} and \ref{Cor:eta} and using
Corollary \ref{Cor1.4}, one may compute the first derivative of an expression
of the form
$$\int_{t\in\Delta_n}M\LP\exp(t_0x),\exp(t_1x),\ldots,\exp(t_nx)\RP
\,d\lambda_n(t),$$
where $M:A^{n+1}\to A$ is bounded multilinear. Its derivative, computed
at $v\in A$ will be
$$\sum_{k=1}^{n+1}\int_{t\in\Delta_{n+1}}M\LP\exp(t_0x),\ldots,
\exp(t_{k-2}x),\exp(t_{k-1}x)v\exp(t_kx),\exp(t_{k+1}x),\ldots,
\exp(t_{n+1}x)\RP\,d\lambda_{n+1}(t).$$
Since this is a similar expression, one may use the same formula to
differentiate again. In this way one gets (by induction on $n$) the formula
for the $n$-th derivative of $\exp x$, namely:

\begin{Thm}
\label{Th3.4} Let $A$ be a \Bh\ algebra. For $x\in A$,
$v_1,\ldots,v_n\in A$,
\begin{eqnarray}
&&\LA D^n(\exp x),v_1\otimes\ldots\otimes v_n\RA=\nonumber\\
&&\sum_{\phi\in S_n}\int_{\Delta_n}
\exp(t_0x)v_{\phi(1)}\exp(t_1x)v_{\phi(2)}\ldots\exp(t_{n-1}x)v_{\phi(n)}
\exp(t_nx)\,d\lambda_n(t).\label{(3.7)}
\end{eqnarray}
\end{Thm}\qed

\section[Applying a $\bC^\I$-\fn\ to Hermitians in a $C^*$-algebra]%
{Applying a $\bC^\I$-\fn\ to Hermitian elements of a $C^*$-algebra}
\label{Sec4}
In this section we assume $A$ is a unital $C^*$-algebra. Denote by $A_0$
the real \Bh\ space of Hermitian elements in $ A$. For $x\in A_0$ denote
by $B_x$ the unital $C^*$-algebra generated by $x$, which is commutative,
and by $\spc x$ the spectrum of $x$ w.r.t.\ $B_x$. $\spc x$ is a compact
subset of $\bR$ and $ B_x$ is isomorphic as a $C^*$-algebra to the
$C^*$-algebra $\cC(\spc x)$ of all \cs\ complex \fn s on $\spc x$, indeed
there is a unique $*$-isomorphism $\Phi_x:\cC(\spc x)\to B_x$ mapping the
\fn\ $t$ (the inclusion \fn\ to $\bC$) to $x$.
\parsk

{\EM Remark} $\spc x$ is also the spectrum of $x$ w.r.t.\ $A$. Indeed,
for $t_0\in\bC$, if $t_0\notin\spc x$ then $x-t_0\cdot1$ is invertible in
$B_x$ hence in $A$. If $t_0\in\spc x$ then there are $g_n\in\cC(\spc x)$
so that $\sup_{t\in\spc x}|g_n(t)|=1$ but 
$\sup_{t\in\spc x}|(t-t_0)g_n(t)|\to_n0$. Hence $\|\Phi_x(g_n)\|=1$
and $\|(x-t_0\cdot1)\Phi_x(g_n)\|\to_n0$. If $x-t_0\cdot1$ were invertible
in $ A$ that would imply
$$1=\|\Phi_x(g_n)\|\le\|(x-t_0\cdot1)^{-1}\|
\|(x-t_0\cdot1)\Phi_x(g_n)\|\to_n0.$$
\parsk

Let $g:\bR\to\bC$ be \cs. For $x\in A_0$ we denote by $g$ also the restriction
of $g$ to $\spc x$, and define $\Phi_x(g)\in A$ as {\EM the value of $g$
applied to $x$} and denote it by $g(x)$. Thus $g$ induces a \fn\ $A_0\to A$,
which will be denoted by $g_*$.
\parsk

\begin{Prop}
\label{Prop4.1} Let $x\in A_0$ and $s\in\bR$. Then the \fn\
$t\mapsto\exp(ist)$ applied to $x$ gives $\exp(isx)$ in the sense of
\S\ref{Sec3}, and $\|\exp(isx)\|=1$.
\end{Prop}

\begin{Prf}
The complex exponential series converges uniformly on compacta. Therefore,
uniformly for $t\in\spc x$, we have
$$\exp(ist)=\sum_{k=0}^\I\dfrac1{k!}(is)^kt^k.$$
Applying the \cs\ homomorphism $\Phi_x$ which maps the \fn\ $t$ to $x$,
one obtains
$$\Phi_x(t\mapsto\exp(ist))=\sum_{k=0}^\I\dfrac1{k!}(is)^kx^k,$$
which is $\exp(isx)$ in the sense of \S\ref{Sec3}. Also, 
$\spc x\subset\bR$ hence $\exp(ist)$ has for each $t$ the absolute value $1$,
so its norm is $1$ and since $\Phi_x$ is an isometry, $\|\exp(isx)\|=1$.
\end{Prf}
\parsk

We shall use the {\EM Fourier transform} on $\bR$. For $g\in L^1(\bR)$ we
take its Fourier transform as
\BE\HAT{g}(s):=\dfrac1{2\pi}\int_\bR\exp(-ist)g(t)\,dt.\label{(4.1)}\EE
(This means that we take on $\bR$ the Haar measure which is $1/(2\pi)$
times Lebesgue, so on the dual $s$-line we will have to take the Haar
measure = Lebesgue.)
It is well-known that for $g\in L^1(\bR)$,\,\,$\HAT{g}\in\cC_0(\bR)$, and
if $\HAT{g}$ also happens to be in $L^1(\bR)$ then
\BE g(t)=\int_\bR\exp(ist)\HAT{g}(s)\,ds.\label{(4.2)}\EE
\parsk

Consider Schwartz's space $\cS$ -- the linear space of complex $\cC^\I$
\fn s $g(t)$ on $\bR$ such that for all nonnegative integers $n$, $k$,
$|t|^k|g^{(n)}(t)|$ is bounded, with Fr\'echet space structure defined
by the seminorms $\sup_t|t|^k|g^{(n)}(t)|$ for all $n$, $k$. It is well-known
that if $g\in\cS$ then $\HAT{g}\in\cS$ and $(\ref{(4.2)})$ holds, and that
the mapping $g\mapsto\HAT{g}$ is a topological automorphism of $\cS$.
\parsk

Denote $\cC^\I:=\cB\cC^\I(\bR,\bC)$. This is the Fr\'echet space of all
complex (generally unbounded) $\cC^\I$ \fn s on $\bR$. Denote by
$\LP\cC^\I\RP^*$ the dual space, whose members may be viewed as the 
complex distributions on $\bR$ with compact support.
\parsk

\begin{Thm}
\label{Th4.2} If $g\in\cC^\I$ then $g_*:A_0\to A$ belongs to 
$\cB\cC^\I(A_0,A)$.
\end{Thm}

\begin{Prf}
We have to prove $g_*$ is $\cC^\I$ with bounded derivatives on $\|x\|<r$ for
all $r$. But for $\|x\|<r$, $g(x)=g_*(x)$ depends only on the restriction
of $g$ to $\spc x$ contained in $[-r;r]$ hence it is equal to $g_1(x)$ if
$g_1\in\cC^\I$ coincides with $g$ on $[-r;r]$. We may take $g_1$ with
compact support. Hence it suffices to prove the proposition for $g$ with
compact support. Such $g$ is in $\cS$, has a Fourier transform $\HAT{g}$
and (\ref{(4.2)}) holds pointwise in $t$.
\parsk

Note, that by \S\ref{Sec2} and \S\ref{Sec3}, for each $s\in\bR$
$x\mapsto\exp(isx)$ belongs to $\cB\cC^\I(A_0,A)$ and
$$\dfrac{d^n}{dx^n}\exp(isx)=(is)^nD^n(\exp y)|_{y=isx}.$$
For $x\in A_0$, Thm.\ \ref{Th3.4} gives a formula (\ref{(3.7)}) for
$D^n(\exp y)$ and Prop.\ \ref{Prop4.1} says that the exponents appearing
in the integral in (\ref{(3.7)}) have norm $1$. Also, $S_n$ has $n!$ elements
and $\lambda_n(\Delta_n)=1/(n!)$. Thus one concludes that
$\|D^n(\exp y)|_{y=isx}\|\le1$ hence
$$\left\|\dfrac{d^n}{dx^n}\exp(isx)\right\|\le|s|^n.$$
So, if we denote the element $x\mapsto\exp(isx)$ of $\cB\cC^\I(A_0,A)$ by
$e_s$, then $\rho_{n,r}(e_s)\le|s|^n$ independently of $r$. Also, by Remark
\ref{Rmk2.4} the \fn\ $s\mapsto e_s$ is \cs\ from $\bR\to\cB\cC^\I(A_0,A)$.
Return now to $g$. Since $\HAT{g}\in\cS$, $\int|s|^n|\HAT{g}(s)|\,ds<\I$
i.e.\ $\int\rho_{n,r}(e_s)|\HAT{g}(s)|\,ds<\I$. Therefore in
$\cB\cC^\I(A_0,A)$ the following integral exists:
\BE g_1=\int_\bR e_s\HAT{g}(s)\,ds\label{(4.3)}\EE
and, of course, $g_1\in\cB\cC^\I(A_0,A)$. For any $x\in A_0$ we may apply
to (\ref{(4.3)}) the evaluation operator mapping each element of 
$\cB\cC^\I(A_0,A)$ to its value at $x$, to obtain
\BE g_1(x)=\int_\bR\exp(isx)\HAT{g}(s)\,ds.\label{(4.4)}\EE
On the other hand, for $s\in\bR$ the element of $\cC(\spc x)$
$t\mapsto\exp(ist)$ has norm $1$ and by Prop.\ \ref{Prop1.2} the \fn\ mapping
each $s$ to this element is \cs. Therefore the integral
$\exp(ist)\HAT{g}(s)\,ds$ exists in $\cC(\spc x)$ and by (\ref{(4.2)}) has
for each $t\in\spc x$ the value $g(t)$. Hence the value of this integral
in $\cC(\spc x)$ is (the restriction of) $g$, and applying $\Phi_x$ gives
in $B_x$ hence in $A$
$$g(x)=\int_\bR\exp(isx)\HAT{g}(s)\,ds$$
which with (\ref{(4.4)}) gives $g(x)=g_1(x)$ hence $g_*=g_1$ implying
$g_*\in\cB\cC^\I(A_0,A)$.
\end{Prf}
\parsk

\section{Bounding seminorms of $g_*$ by the seminorms of $g$}
\label{Sec5}
We first consider $\cC^\I$ and its dual $\LP\cC^\I\RP^*$ --
the space of complex distributions with compact support.
\parsk

For $G\in\LP\cC^\I\RP^*$, define its {\EM Fourier transform} as the \fn\
on $s\in\bC$
\BE\HAT{G}(s):=\LA G,t\mapsto\exp(-ist)\RA.\label{(4.5)}\EE
\parsk

By the definition of the topology in $\cC^\I$ and Prop.\ \ref{Prop1.3},
the \fn\ mapping $s\in\bC$ to $[t\mapsto\exp(-ist)]\in\cC^\I$ is 
differentiable in $s$ with derivative the operator ($v\in\bC$)
$v\mapsto[t\mapsto-it\exp(-ist)v]$. (Apply Prop.\ \ref{Prop1.3} separately for
derivatives of each order and for each compact interval of $t$.) This means that
this mapping, as a $\cC^\I$-valued \fn\ on $\bC$, is entire. Consequently
also $\HAT{G}$ is entire. Moreover, one easily finds that in $\cC^\I$, for
{\em real} $s$, $\rho_{n,k}(t\mapsto\exp(-ist))\le|s|^k$, hence for $s\in\bR$
\,$\HAT{G}$ is bounded by a polynomial. 
\parsk

Suppose now that $g\in\cS$. Then since $\int|s|^k|\HAT{g}(s)|\,ds<\I$,
one finds that $\int\HAT{g}(s)\exp(-ist)\,ds$ exists in $\cC^\I$ and its
value in each $t$ is, by $(\ref{(4.2)})$, $g(-t)$, hence its value in
$\cC^\I$ must be the \fn\ $t\mapsto g(-t)$. Applying $G$ one obtains:
$$\int_\bR\HAT{g}(s)\HAT{G}(s)\,ds=\LA G,t\mapsto g(-t)\RA\qquad
G\in\LP\cC^\I\RP^*,g\in\cS,$$
or, written differently,
\BE\LA G,g\RA=\int_\bR\HAT{G}(s)\HAT{g}(-s)=
\int_\bR\HAT{G}(-s)\HAT{g}(s).
\qquad G\in\LP\cC^\I\RP^*,g\in\cS.\label{(gG)}\EE
\parsk

Let now $G\in\LP\cC^\I\RP^*$, $g\in\cC^\I$. The linear operator $f\mapsto fg$
from $\cC^\I$ to itself is \cs, and we define (as usual) the product
$gG\in\LP\cC^\I\RP^*$ by
$$\LA gG,f\RA:=\LA G,gf\RA\qquad f\in\cC^\I.$$
For $g\in\cS$, we have
$$\HAT{(gG)}(u)=\LA gG,t\mapsto\exp(-iut)\RA=
\LA G,t\mapsto g(t)\exp(-iut)\RA.$$
We wish to use $(\ref{(gG)})$. The Fourier transform of 
$t\mapsto g(t)\exp(-iut)$ is
$$s\mapsto\dfrac1{2\pi}\int_\bR\exp(-ist)g(t)\exp(-iut)\,dt=\HAT{g}(u+s)$$
and using $(\ref{(gG)})$ we find
\BE\HAT{(gG)}(u)=\int_\bR\HAT{G}(s)\HAT{g}(u-s)\,ds\qquad
G\in\LP\cC^\I\RP^*,g\in\cS\label{(4.10)}\EE
i.e.\ $\HAT{(gG)}$ is the convolution of $\HAT{G}$ and $\HAT{g}$.
\parsk

We turn now to derivatives and integrals.
Denote by $\bf1$ the element of $\cC^\I$ which is the constant \fn\
$1$ and by $\delta$ the element of $\LP\cC^\I\RP^*$ given by ($g\in\cC^\I$)
$\LA\delta,g\RA:=g(0)$.
\parsk

The linear operator mapping each $g\in\cC^\I$ to its derivative is \cs.
Also the linear operator mapping $g$ to
$\int_0g:=\LQ t\mapsto\int_0^tg(\tau)\,d\tau\RQ$
($t$ may also be zero or negative) is \cs. For $G\in\LP\cC^\I\RP^*$, define
$G'$ (as usual) and $\int_0G$ by
\BE\LA G',g\RA:=-\LA G,g'\RA\qquad\LA\int_0G,g\RA:=-\LA G,\int_0g\RA\qquad
g\in\cC^\I\label{(4.6)}\EE
Note that if $\supp G$ (the support of $G$) is contained in $[-r;r]$, also
$\supp G'$ and $\supp\int_0G$ are contained in $[-r;r]$.
\parsk

We have
$$\LA\LP\int_0G\RP',g\RA=
-\LA\int_0G,g'\RA=
\LA G,\int_0(g')\RA=
\LA G,g-\LA\delta,g\RA{\bf1}\RA=
\LA G,g\RA-\LA G,{\bf1}\RA\LA\delta,g\RA.$$
Hence
\BE\LP\int_0G\RP'=G-\LA G,{\bf1}\RA\delta.\label{(4.7)}\EE

For the Fourier transforms we have, by $(\ref{(4.5)})$ and $(\ref{(4.6)})$
$$\HAT{(G')}(s)=\LA G',t\mapsto\exp(-ist)\RA=
-\LA G,t\mapsto-is\exp(-ist)\RA=is\LA G,t\mapsto\exp(-ist)\RA=is\HAT{G}(s).$$
Thus
\BE\HAT{(G')}(s)=is\HAT{G}(s).\label{(4.8)}\EE
And by $(\ref{(4.7)})$,
$$is\HAT{\LP\int_0G\RP}(s)=G(s)-\LA G,{\bf1}\RA\HAT{\delta}(s),$$
and since $\HAT\delta\equiv1$ and $\LA G,{\bf1}\RA=\HAT{G}(0)$, one concludes
that for $s\ne0$
\BE\HAT{\LP\int_0G\RP}(s)=
\dfrac1{is}\LP\HAT{G}(s)-\HAT{G}(0)\RP.\label{(4.9)}\EE
\par\bigskip

Now let $g\in\cC^\I$ $r>0$,\,\,$n=0,1,2,\ldots$ and we wish to bound the
seminorm in $\cB\cC^\I(A_0,A)$:\,\,\,$\rho_{n,r}(g_*)$. (The integrals in
what follows are always on the whole $\bR$.)
\parsk

Let $f\in\cS$ with $f\equiv g$ on $[-r;r]$. For $x\in A_0$ with
$\|x\|\le r$ we have $g(x)=f(x)$, Therefore $\rho_{n,r}(g_*)=\rho_{n,r}(f_*)$.
With the notation of the proof of Thm.\ \ref{Th4.2}, we have,
in $\cB\cC^\I(A_0,A)$, $f_*=\int e_s\HAT{f}(s)\,ds$ hence
$$\rho_{n,r}(f_*)\le\int\rho_{n,r}(e_s)\left|\HAT{f}(s)\right|\,ds\le
\int|s|^k\left|\HAT{f}(s)\right|\,ds.$$
Thus, if we denote by $Y$ the closed affine subspace of $\cS$:
$$Y:=\LB f\in\cS\,\big|\,f\equiv g\mbox{ on }[-r;r]\RB,$$
then
\BE\rho_{n,r}(g_*)\le\inf_{f\in Y}\int|s|^n\left|\HAT{f}(s)\right|\,ds
\label{(4.16)}\EE
Denote the infimum in the RHS by $a$. If $a=0$ then $\rho_{n,r}(g_*)=0$.
If $a>0$ then $Y$ is disjoint from the set
$$U:=\LB h\in\cS\,\big|\,\int|s|^n\left|\HAT{f}(s)\right|\,ds<a\RB.$$
$U$ is convex, balanced and open in $\cS$ (since
$\int|s|^n\left|\HAT{f}(s)\right|\,ds$ is a \cs\ seminorm there).
By a form of Hahn-\Bh\ $\exists$ a (closed) hyperplane in $\cS$ containing $Y$
and disjoint from $U$. In other words, $\exists$ a \cs\ linear \fn al $G$ on
$\cS$ (i.e.\ a temperate distribution) so that $G<a$ on $U$ and $G=a$ on $Y$.
In particular, $G$ vanishes on differences of members of $Y$, i.e.\ on
members of $\cS$ vanishing identically on $[-r;r]$, which means that
$\supp G\subset[-r;r]$. Thus $G$ may be seen as a member of $\LP\cC^\I\RP^*$.
\parsk

Using $(\ref{(gG)})$ one may express $G$ being $<a$ on $U$ by:
$$\forall h\in\cS\quad
\int|s|^n\left|\HAT{h}(s)\right|\,ds<a\Rightarrow
\left|\int\HAT{G}(-s)\HAT{h}(s)\,ds\right|<a$$
and since both sides are linear in $h$, we find
\BE\forall h\in\cS\quad
\left|\int\HAT{G}(-s)\HAT{h}(s)\,ds\right|\le
\int|s|^n\left|\HAT{h}(s)\right|\,ds.\label{(4.17)}\EE
This implies $\left|\HAT{G}(s)\right|\le|s|^n$. If $n\ge1$ then $\HAT{G}(0)=0$
and by $(\ref{(4.9)})$,\,$\left|\HAT{\LP\int_0G\RP}(s)\right|\le|s|^{n-1}$
for $s\ne0$ and by continuity of $\HAT{\LP\int_0G\RP}$ also for $s=0$.
Since $\LA G,{\bf1}\RA=\HAT{G}(0)=0$, one find by $(\ref{(4.7)})$ that
$G=\LP\int_0G\RP'$. Also $\supp\int_0G\subset[-r;r]$. Thus we found that if
$n\ge1$ then $G$ is a derivative of a member of $\LP\cC^\I\RP^*$ with
support contained in $[-r;r]$ and Fourier transform bounded by $|s|^{n-1}$
on $\bR$.
Proceeding by induction, we find that $G$ is the $n$-th
derivative of a $G_1\in\LP\cC^\I\RP^*$ with $\supp G_1\subset[-r;r]$
and $\left|\HAT{G_1}(s)\right|\le1$ for real $s$.
\parsk

Choose an $f\in Y$, i.e.\ $f\equiv g$ on $[-r;r]$. We have $|\LA G,f\RA|=a$
hence $|\LA G,g\RA|=a$, i.e.\
$a=\left|\LA G_1^{(n)},g\RA\right|=\left|\LA G_1,g^{(n)}\RA\right|$. Recall
that by $(\ref{(4.16)})$ $a$ was an upper bound to $\rho_{n,r}(g_*)$. Thus we
conclude:

\begin{Thm} \label{Th:rhoG}
Let $g\in\cC^\I$, $r>0$ and $n=0,1,\ldots$. Then $\exists$ a
$G\in\LP\cC^\I\RP^*$ with $\supp G\subset[-r;r]$ and
$\left|\HAT{G}(s)\right|\le1$ for $s\in\bR$ so that
$$\rho_{n,r}(g_*)\le\left|\LA G,g^{(n)}\RA\right|.$$
\end{Thm}\qed
\parsk

To obtain an explicit bound, we apply:

\begin{Main}
{\EM Bernshtein's Theorem}: If $G$ is as in Thm.\ \ref{Th:rhoG} then for
real $s$ $\left|\LP\HAT{G}\RP'(s)\right|\le r$.
\end{Main}
\begin{Prff} {\EM Proof of Bernshtein's Theorem}:
We have seen that the derivative of the \fn\ $\bR\to\cC^\I$ mapping $s$
to $t\mapsto\exp(-ist)$ is the \fn\ mapping $s$ to $t\mapsto-it\exp(-ist)$,
hence
\BE\LP\HAT{G}\RP'(s)=\LA G,t\mapsto-it\exp(-it)\RA.\label{(5.1)}\EE
Consider the circle $\bT=\bR/(2\pi\bZ)$. We use the same notation for a
$t\in\bR$ and its image in $\bT$. In $\bT$ we take a Haar measure normalized
to mass $1$ for the whole $\bT$. Consider the \fn\ on $\bT$:
\BE j(t):=\left\{
\begin{array}{ll}\dfrac\pi2-t&0\le t\le\pi\\ \\
\dfrac\pi2+t&-\pi\le t\le0\end{array}\right.\label{(5.2)}\EE
$j$ is even, \cs, and is easily seen to be $-\pi/2$ plus $2\pi$ times the
convolution of the characteristic \fn\ of the interval $[-\pi/2;\pi/2]$
with itself. Hence all Fourier coefficients $\HAT{j}(m)$,\,$m\ne0$ are
non-negative, while clearly $\HAT{j}(0)=0$.
\parsk

Let $\pi/2>\eps>0$. Choose a \fn\ $f_\eps$ on $\bT$ non-negative, even,
$\cC^\I$, supported in $[-\eps;\eps]$, with integral $1$ and with
non-negative Fourier coefficients (the last requirement can be acheived
by replacing $f_\eps$ with $f_{\eps/2}$ convolved with itself). Certainly,
its Fourier coefficients $\HAT{\LP f_\eps\RP}(m)$ are fastly decaying,
i.e.\ decay faster than any negative power of $m$.
\parsk

Let $h$ be the convolution $f_\eps*j$. $h$ is real, even, \cs, and
$\HAT{h}(m)$ are non-negative and fastly decaying, certainly in $\ell^1$,
hence
\BE h(t)=\sum_m\HAT{h}(m)\exp(imt)\label{(5.3)}\EE
Moreover, for $\eps<t<\pi-\eps$ $f_\eps$ ``sees'' only the part 
$\dfrac\pi2-t$ of $j$, and since $f_\eps$ is even $h(t)=\dfrac\pi2-t$.
\parsk

Let $c<\pi/(2r)$ and subsitute $ct+\pi/2$ for $t$ in $(\ref{(5.3)})$:
\BE h(ct+\pi/2)=\sum_m\HAT{h}(m)\exp(im(ct+\pi/2)),\label{(5.6)}\EE
and view this as a series of members of $\cC^\I$ (\fn s on $\bR$).
On the one hand, any seminorm on $\exp(im(ct+\pi/2))$ grows polynomially in
$m$ and the coefficients are fastly decaying, hence the series converges
in $\cC^\I$ to $h(ct+\pi/2)$. On the other hand, for $\eps$ small enough
$t\mapsto ct+\pi/2$ will project $[-r;r]$ into the interval where
$h(t)=\dfrac\pi2-t$ hence $h(ct+\pi/2)=-ct$. This holds on a neighborhood
of $\supp G$. Multiplying by $\dfrac{i}{c}\exp(-ist)$ we have a series
converging in $\cC^\I$, such that in a neighborhood of $\supp G$:
$$-it\exp(-ist)=\dfrac{i}{c}\sum_m\HAT{h}(m)\exp(i(mc-s)t)\exp(im\pi/2)$$
Applying $G$, we get in the LHS $\LP\HAT{G}\RP'(s)$ and in the RHS a
combination
of $\exp(im\pi/2)\HAT{G}(mc-s)$, which have absolute value $\le 1$ simce
$|\HAT{G}|\le1$ on $\bR$, with the non-negative coefficients $\HAT{h}(m)$.
Therefore
$$\left|\LP\HAT{G}\RP'(s)\right|\le
(1/c)\sum_m\HAT{h}(m)=(1/c)h(0)=(1/c)(\pi/2).$$
This holds for any $c<\pi/2r$, and one concludes that
$\left|\LP\HAT{G}\RP'(s)\right|\le r$.
\end{Prff}
\parsk

Thus for real $s$,
for $G$ as in Thm.\ \ref{Th:rhoG} $\left|\LP\HAT{G}\RP'(s)\right|\le r$,
hence
$$\dfrac1{|s|}\left|\HAT{G}(s)-\HAT{G}(0)\right|\le r$$
while on the other hand, since $|\HAT{G}|\le1$,
$$\dfrac1{|s|}\left|\HAT{G}(s)-\HAT{G}(0)\right|\le\dfrac2{|s|}.$$
Thus $(-i/s)\LP\HAT{G}(s)-\HAT{G}(0)\RP$,
which By $(\ref{(4.9)})$ is just $\HAT{\LP\int_0G\RP}$, is an $L^2$
\fn\ on $\bR$ with square $L^2$ norm (w.r.t.\ Lebesgue) no more than
$$2\int_0^{2/r}r^2\,ds+2\int_{2/r}^\I\LP2/s\RP^2\,ds=
2\cdot\dfrac2r\cdot r^2+2\cdot4\cdot(2/r)^{-1}=8r.$$
We need to estimate $\LA G,g\RA$ for $g\in\cC^\I$. To this end choose
$r_2>r_1>r$ and a $\cC^\I$ \fn\ $\ell$ on $\bR$ such that $0\le\ell(t)\le1$,
$\ell\equiv1$ on $[-r_1;r_1]$ and $\ell\equiv0$ outside $]-r_2;r_2[$.
By $(\ref{(4.7)})$
$$\left|\LA G,g\RA\right|=
\left|\LA G,{\bf1}\RA\LA\delta,g\RA+\LA\LP\int_0G\RP',g\RA\right|\le
|\HAT{G}(0)||g(0)|+\left|\LA\int_0G,g'\RA\right|.$$
$g'$ and $\ell\cdot g'$ coincide on $\supp\int_0G$. By $(\ref{(gG)})$
$$\left|\LA\int_0G,g'\RA\right|=\left|\LA\int_0G,\ell\cdot g'\RA\right|=
\left|\int_\bR\HAT{\LP\int_0G\RP}(-s)\HAT{\LP\ell\cdot g'\RP}(s)\,ds\right|\le
\left\|\HAT{\LP\int_0G\RP}\right\|_2
\left\|\HAT{\LP\ell\cdot g'\RP}\right\|_2,$$
and we know that $\left\|\HAT{\LP\int_0G\RP}\right\|_2\le\sqrt{8r}$, while
by Plancherel's Theorem $\left\|\HAT{\LP\ell\cdot g'\RP}\right\|_2$ is
the $L^2$ norm of $\ell(t)\cdot g'(t)$ (computed w.r.t.\ $1/(2\pi)$ times
Lebesgue) which with $r_1$ and $r_2$ tending to $r$ gives 
$$\left|\LA G,g\RA\right|\le|g(0)|+\sqrt{8r}
\LP\dfrac1{2\pi}\int_{-r}^r|g'(t)|^2\RP^{1/2}.$$
Combining everything with Thm.\ \ref{Th:rhoG} one finally has

\begin{Thm}
Let $g\in\cC^\I$, $r>0$ and $n=0,1,\ldots$. Then
$$\rho_{r,n}(g_*)\LP=\sup_{\|x\|<r}\|g_*^{(n)}(x)\|\RP
\le\left|g^{(n)}(0)\right|+\sqrt{8r}
\LP\dfrac1{2\pi}\int_{-r}^r\left|g^{(n+1)}(t)\right|^2\RP^{1/2}.$$
Clearly, the same will hold also for $g$ in the Sobolev space which is the
completion of the $\cC^\I$-\fn s w.r.t.\ the RHS norm.
\end{Thm}\qed

\end{document}